\documentclass[12pt]{article}
\usepackage{mathrsfs}
\usepackage{amssymb} \textwidth 150mm \textheight 220mm \oddsidemargin
15pt \evensidemargin 0pt \topmargin 0cm \headsep 0.3cm

\usepackage{amsmath}
\usepackage{graphicx}
\usepackage{indentfirst}
\usepackage{cite}
\usepackage{float}
\usepackage{CJK}

\begin{document}
\title{\Large\bf{Bifurcation of limit cycles of the nongeneric quadratic reversible system with discontinuous perturbations
\thanks{E-mail address: jihua1113@163.com(J. Yang)} }}
\author{{Jihua Yang}
\\ {\small \it School of Mathematics and Computer Science, Ningxia Normal University, Guyuan}\\
 {\small \it   756000, China}}
\date{}
\maketitle \baselineskip=0.9\normalbaselineskip \vspace{-3pt}
\noindent
{\bf Abstract}\, By using the Picard-Fuchs equation and the property of Chebyshev space to the discontinuous differential system, we obtain an upper bound of the number of limit cycles for the nongeneric quadratic reversible system when it is perturbed inside all discontinuous polynomials with degree $n$.
\vskip 0.2 true cm
\noindent
\vskip 0.2 true cm


 \section{Introduction and main results}
 \setcounter{equation}{0}
\renewcommand\theequation{1.\arabic{equation}}

Stimulated by discontinuous phenomena in the real world, such as biology \cite{K}, nonlinear oscillations \cite{T}, impact and friction mechanics \cite{BBC}, a big interest has appeared for studying the number of limit cycles and their relative positions of discontinuous differential systems. Similar to the smooth differential system, one of the main problems in the qualitative theory of non-smooth differential systems is the study of their limit cycles, and many methodologies have been developed, such as Abelian integral method (or first order Melnikov function) \cite{LHR,LH,X,YZ16}, averaging method \cite{CLZ,DL,LL,LM,LMN,LT}. This problem can be seen as an extension of the infinitesimal Hilbert's 16th Problem  to the discontinuous world.

The list of quadratic center at (0,0), almost all the orbits of which are cubic, looks as follows \cite{I,ZLLZ}:
 \vskip 0.2 true cm
The Hamiltonian system $Q_3^H$: $\dot{z}=-iz-z^2+2|z|^2+(b+ic)\bar{z}^2$.
\vskip 0.2 true cm

The Hamiltonian triangle: $\dot{z}=-iz+\bar{z}^2$.
\vskip 0.2 true cm

The reversible system: $\dot{z}=-iz+(2b+1)z^2+2|z|^2+b\bar{z}^2$,\ $b\neq-1$.
 \vskip 0.2 true cm

The generic Lotka-Volterra system: $\dot{z}=-iz+(1-ci)z^2+ci\bar{z}^2$,\ $c=\pm\frac{1}{\sqrt{3}}$.
\vskip 0.2 true cm

 Under the perturbations of continuous polynomials of degree $n$, Horozov and Iliev \cite{HI} proved that the number of limit cycles for
$Q_3^H$ and Hamiltonian triangle does not exceed $5n+15$, and Zhao et al. \cite{ZLLZ} proved that the number of limit cycles for
reversible and generic Lotka-Volterra systems does not exceed $7n$.

Let $z=x+iy$ and by a linear transformation, the reversible system can be written \cite{ZLLZ}:
\begin{eqnarray}
\begin{cases}
\dot{x}=xy,\\
\dot{y}=\frac{3}{2}y^2+ax^2-2(a+1)x+a+2,
\end{cases}
\end{eqnarray}
where $a\in\mathbb{R}$. When $a=-2$, system (1.1) corresponds to the nongeneric case of the reversible system (1.1):
\begin{eqnarray}
\begin{cases}
\dot{x}=xy,\\
\dot{y}=\frac{3}{2}y^2-2x^2+2x
\end{cases}
\end{eqnarray}
 whose first integral is
\begin{eqnarray}
 H(x,y)=x^{-3}\Big(\frac{1}{2}y^2-2x^2+x\Big)=h,\ h\in(-1,0)
 \end{eqnarray}
with the integrating factor $\mu(x,y)=x^{-4}$.

In the present paper, by using the Picard-Fuchs equation and the property of Chebyshev space, we investigate the number of limit cycles of system (1.2) under discontinuous polynomial perturbations of degree $n$. The system (1.2) has a center (1,0) and $h=-1$ corresponds to the center (1,0). 
The perturbed system of (1.2) is
\begin{eqnarray}
\left(
  \begin{array}{c}
          \dot{x} \\
          \dot{y}
 \end{array}
 \right)=\begin{cases}
 \left(
  \begin{array}{c}
          xy+\varepsilon f^+(x,y) \\
          \frac{3}{2}y^2-2x^2+2x+\varepsilon g^+(x,y)
 \end{array}
 \right), \quad y>0,\\
 \left(
  \begin{array}{c}
         xy+\varepsilon f^-(x,y) \\
          \frac{3}{2}y^2-2x^2+2x+\varepsilon g^-(x,y)
 \end{array}
 \right),\quad y<0,\\
 \end{cases}
       \end{eqnarray}
 where $0<|\varepsilon|\ll1$,
$$f^\pm(x,y)=\sum\limits_{i+j=0}^na^\pm_{i,j}x^iy^j,\ \ g^\pm(x,y)=\sum\limits_{i+j=0}^nb^\pm_{i,j}x^iy^j,\ i,j\in\mathbb{N}.$$

From the Theorems 1.1 in \cite{LH,HS}, by linear transformations, we know that the first order Melnikov function $M(h)$ of system (1.4) is
\begin{eqnarray}
\begin{aligned}
M(h)=&\int_{\Gamma^+_h}x^{-4}[g^+(x,y)dx-f^+(x,y)dy]\\
&+\int_{\Gamma^-_h}x^{-4}[g^-(x,y)dx-f^-(x,y)dy],\ \  h\in(-1,0),
\end{aligned}
\end{eqnarray}
where
$$\begin{aligned}&\Gamma_h^+=\{(x,y)|H(x,y)=h,h\in(-1,0),y>0\},\\
&\Gamma^-_h=\{(x,y)|H(x,y)=h,h\in(-1,0),y<0\},\end{aligned}$$
and its number of zeros gives an upper bound of the number of limit cycles of system (1.4) bifurcating from the period annulus.

Our main results are the following two theorems.
 \vskip 0.2 true cm

\noindent
{\bf Theorem 1.1.}\, {\it Suppose that $h\in(-1,0)$.}
 \vskip 0.2 true cm
\noindent
(i) {\it If $n=2,3$, then the number of limit cycles of system (1.4) bifurcating from the period annulus is not more than 40 (counting multiplicity). }
 \vskip 0.2 true cm
 \noindent
(ii) {\it If $4\leq n\leq 7$, then the number of limit cycles of system (1.4) bifurcating from the period annulus is not more than $24n-56$ (counting multiplicity).}
 \vskip 0.2 true cm
 \noindent
(iii) {\it If $n\geq 8$,  then the number of limit cycles of system (1.4) bifurcating from the period annulus is not more than $22n-64$ (counting multiplicity).}

 \vskip 0.2 true cm

\noindent
{\bf Theorem 1.2.}\, {\it Suppose that $h\in(-1,0)$, $a^+_{i,j}=a^-_{i,j}$ and $b^+_{i,j}=b^-_{i,j}$.}
 \vskip 0.2 true cm
\noindent
(i) {\it If $n=2,3$,  then the number of limit cycles of system (1.4) bifurcating from the period annulus is not more than $4$ (counting multiplicity).}
 \vskip 0.2 true cm
 \noindent
(ii) {\it If $n\geq4$,  then the number of limit cycles of system (1.4) bifurcating from the period annulus is not more than $3n-8$ (counting multiplicity).}
\vskip 0.2 true cm

\noindent
{\bf Remark 1.1}. (i) By using the Picard-Fuchs equation, we greatly simplified the computation of the first order Melnikov function. And then we can estimate the number of zeros of the first order Melnikov function which controls the number of limit cycles of the corresponding perturbed system benefited from the property of Chebyshev space. It is worth noting that these methods can be applied to study the bifurcation of limit cycles for other integrable differential systems.
\vskip 0.2 true cm

\noindent
(ii) The perturbation as in (1.4) can be found in many practical applications, such as in the slender rocking block model and nonlinear compliant oscillator, see \cite{SD,PLBL,H} and the references quoted there.
\vskip 0.2 true cm

\noindent
(iii) If $h\in(-1,0)$, $a^+_{i,j}=a^-_{i,j}$ and $b^+_{i,j}=b^-_{i,j}$, then Zhao et al. \cite{ZLLZ} obtained that the number of limit cycles of system (1.4) bifurcating from the period annulus is not more than $3n-4$ for $n\geq4$; 8 for $n=3$; 5 for $n=2$ (counting multiplicity).
\vskip 0.2 true cm

The rest of the paper is organized as follows: In Section 2, we will obtain the algebraic structure of the first order Melnikov function $M(h)$ and the Picard-Fuchs equations satisfied by the generators of $M(h)$ are also obtained. Finally, we will prove Theorems 1.1 and 1.2 in Section 3.

\section{The algebraic structure of $M(h)$ and Picard-Fuchs equation}
 \setcounter{equation}{0}
\renewcommand\theequation{2.\arabic{equation}}

In this section, we obtain the algebraic structure of the first Melnikov function $M(h)$.
For $h\in(-1,0)$, we denote
\begin{eqnarray*}
I_{i,j}(h)=\int_{\Gamma^+_h}x^{i-4}y^jdx,\ \ J_{i,j}(h)=\int_{\Gamma^-_h}x^{i-4}y^jdx.
\end{eqnarray*}
 We first prove the following results.

 \vskip 0.2 true cm

\noindent
{\bf Lemma 2.1.}\, {\it Suppose that $h\in(-1,0)$, $i=-1,0,1,\cdots$ and $j=0,1,2,\cdots$.}
\vskip 0.2 true cm

\noindent
(i)\, The following equalities hold:
\begin{eqnarray}
\begin{cases}
I_{-1,1}(h)=\frac{1}{7}[hI_{1,1}(h)+8I_{0,1}(h)],\\
I_{0,0}(h)=\frac{1}{3}[hI_{2,0}(h)+4I_{1,0}(h)],\\
I_{-1,2}(h)=\frac{4}{3}(h+1)I_{2,0}(h),\\
I_{1,0}(h)=I_{2,0}(h),\\
I_{-1,3}(h)=12[I_{1,1}(h)-I_{0,1}(h)].
\end{cases}
\end{eqnarray}
\begin{eqnarray}
\begin{cases}
I_{-1,4}(h)=4[I_{1,2}(h)-I_{0,2}(h)],\\ I_{0,3}(h)=4[I_{2,1}(h)-I_{1,1}(h)],\\
I_{1,2}(h)=\frac{1}{h}[2I_{0,2}(h)-3I_{-1,2}(h)],\\ I_{2,1}(h)=\frac{1}{h}[4I_{1,1}(h)-5I_{0,1}(h)],\\ I_{3,0}(h)=\frac{1}{h}\big[\frac{1}{2}I_{0,2}(h)-2I_{2,0}(h)+I_{1,0}(h)\big].
\end{cases}
\end{eqnarray}
\vskip 0.2 true cm

\noindent
(ii)\, {If $4\leq n\leq7$, then
\begin{eqnarray*}
\begin{cases}
I_{i,2j+1}(h)=\frac{1}{h^{n-3}}[\bar{\alpha}(h)I_{0,1}(h)+\bar{\beta}(h)I_{1,1}(h)],\ \, i+2j+1=n,\\ I_{i,2j}(h)=\frac{1}{h^{n-3}}[\bar{\gamma}(h)I_{2,0}(h)+\bar{\delta}(h)I_{0,2}(h)],\quad\ \ i+2j=n,\\
\end{cases}
\end{eqnarray*}
 where $\bar{\alpha}(h)$, $\bar{\beta}(h)$, $\bar{\gamma}(h)$ and $\bar{\delta}(h)$ are polynomials of $h$ with $\deg \bar{\alpha}(h),\deg \bar{\delta}(h)\leq n-4$ and $\deg \bar{\beta}(h),\deg \bar{\gamma}(h)\leq n-3$. }
\vskip 0.2 true cm

\noindent
 (iii)\, {If $n\geq8$, then
\begin{eqnarray*}
\begin{cases}
I_{i,2j+1}(h)=\frac{1}{h^{n-3}}[\bar{\alpha}(h)I_{0,1}(h)+\bar{\beta}(h)I_{1,1}(h)],\ \, i+2j+1=n,\\ I_{i,2j}(h)=\frac{1}{h^{n-3}}\bar{\gamma}(h)I_{2,0}(h),\qquad\qquad\qquad\qquad \, i+2j=n,\\
\end{cases}
\end{eqnarray*}
 where $\bar{\alpha}(h)$, $\bar{\beta}(h)$ and $\bar{\gamma}(h)$ are polynomials of $h$ with $\deg \bar{\alpha}(h)\leq n-5$ and $\deg \bar{\beta}(h),\deg \bar{\gamma}(h)\leq n-4$.}

  \vskip 0.2 true cm

\noindent
{\bf Proof.}\, Let $D$ be the interior of $\Gamma_{h}^+\cup \overrightarrow{AB}$, see the black line in Fig.\,1. Using the Green's Formula, we have for $j\geq0$
\begin{eqnarray*}
\begin{aligned}
\int_{\Gamma^+_{h}}x^iy^jdy
=&\oint_{\Gamma^+_{h}\cup\overrightarrow{AB}}x^iy^jdy-\int_{\overrightarrow{AB}}x^iy^jdy\\
=&\oint_{\Gamma^+_{h}\cup\overrightarrow{AB}}x^iy^jdy
= -i\iint\limits_Dx^{i-1}y^jdxdy,
\end{aligned}
\end{eqnarray*}
\begin{eqnarray*}
\begin{aligned}
\int_{\Gamma^+_{h}}x^{i-1}y^{j+1}dx=
\oint_{\Gamma^+_{h}\cup\overrightarrow{AB}}x^{i-1}y^{j+1}dx
=(j+1)\iint\limits_Dx^{i-1}y^jdxdy.
\end{aligned}
\end{eqnarray*}
Hence, \begin{eqnarray}\int_{\Gamma^+_{h}}x^iy^jdy=-\frac{i}{j+1}\int_{\Gamma^+_{h}}x^{i-1}y^{j+1}dx,\ j\geq0.\end{eqnarray}
In a similar way, we have
\begin{eqnarray}\int_{\Gamma^-_{h}}x^iy^jdy=-\frac{i}{j+1}\int_{\Gamma^-_{h}}x^{i-1}y^{j+1}dx,\ j\geq0.\end{eqnarray}
By a straightforward calculation and noting that (2.3) and (2.4), we obtain
\begin{eqnarray}
\begin{aligned}
M(h)=&\int_{\Gamma^+_h}x^{-4}\big(g^+(x,y)dx-f^+(x,y)dy\big)\\
&+\int_{\Gamma^-_h}x^{-4}\big(g^-(x,y)dx-f^-(x,y)dy)\\
=&\int_{\Gamma_h^+}\sum\limits_{i+j=0}^nb^+_{i,j}x^{i-4}y^jdx-\int_{\Gamma_h^+}
\sum\limits_{i+j=0}^na^+_{i,j}x^{i-4}y^jdy\\
&+\int_{\Gamma_h^-}\sum\limits_{i+j=0}^nb^-_{i,j}x^{i-4}y^jdx-\int_{\Gamma_h^-}
\sum\limits_{i+j=0}^na^-_{i,j}x^{i-4}y^jdy\\
=&\sum\limits_{i+j=0}^nb^+_{i,j}\int_{\Gamma^+_h}x^{i-4}y^jdx+
\sum\limits_{i+j=0}^n\frac{i-4}{j+1}a^+_{i,j}\int_{\Gamma^+_h}x^{i-5}y^{j+1}dx\\
&+\sum\limits_{i+j=0}^nb^-_{i,j}\int_{\Gamma^-_h}x^{i-4}y^jdx+
\sum\limits_{i+j=0}^n\frac{i-4}{j+1}a^-_{i,j}\int_{\Gamma^-_h}x^{i-5}y^{j+1}dx\\
=&\sum\limits_{i+j=0,i\geq-1,j\geq0}^n\tilde{a}_{i,j}I_{i,j}(h)
+\sum\limits_{i+j=0,i\geq-1,j\geq0}^n\tilde{b}_{i,j}J_{i,j}(h)\\
:=&\sum\limits_{i+j=0,i\geq-1,j\geq0}^n\rho_{i,j}I_{i,j}(h),
\end{aligned}
\end{eqnarray}
where in the last equality we have used that $J_{i,j}(h)=(-1)^{j+1}I_{i,j}(h)$.

Differentiating (1.3) with respect to $x$, we obtain
\begin{eqnarray}
x^{-3}y\frac{\partial y}{\partial x}-\frac{3}{2}x^{-4}y^2+2x^{-2}-2x^{-3}=0.
\end{eqnarray}
Multiplying (2.6) by $x^iy^{j-2}dx$, integrating over $\Gamma^+_h$ and noting that (2.3), we have
\begin{eqnarray}
(2i+3j-6)I_{i,j}=4j(I_{i+2,j-2}-I_{i+1,j-2}).
\end{eqnarray}
Similarly, multiplying (1.3) by $x^{i-4}y^jdx$ and integrating over $\Gamma^+_h$ yields
\begin{eqnarray}
hI_{i,j}=\frac{1}{2}I_{i-3,j+2}-2I_{i-1,j}+I_{i-2,j}.
\end{eqnarray}
Eliminating $I_{i-3,j+2}$ by (2.7) and (2.8) gives
\begin{eqnarray}
(2i+3j-6)hI_{i,j}=(2i+j-10)I_{i-2,j}-4(i+j-4)I_{i-1,j}.
\end{eqnarray}

From (2.7) we have
\begin{eqnarray}
I_{1,0}=I_{2,0},\ \ I_{-1,3}=12(I_{1,1}-I_{0,1}).
\end{eqnarray}
From (2.8) we obtain
\begin{eqnarray}
hI_{2,0}=\frac{1}{2}I_{-1,2}-2I_{1,0}+I_{0,0}.
\end{eqnarray}
Taking $(i,j)=(2,0),(1,1)$ in (2.9) we have
\begin{eqnarray}
I_{0,0}=\frac{1}{3}(hI_{2,0}+4I_{1,0}),\ I_{-1,1}=\frac{1}{7}(hI_{1,1}+8I_{0,1}).
\end{eqnarray}
Hence,
\begin{eqnarray}
I_{0,0}=\frac{1}{3}(h+4)I_{2,0}.
\end{eqnarray}
From (2.10)-(2.12) we get
\begin{eqnarray}
I_{-1,2}=\frac{4}{3}(h+1)I_{2,0}.
\end{eqnarray}
(2.10) and (2.12)-(2.14) imply (2.1) holds. In a similar way, applying the equalities (2.7) and (2.9), we can obtain
(2.2). Hence, the conclusion (i) holds. By some straightforward calculations according to (2.7) and (2.9), we can get the conclusion (ii).

(iii) Now we prove the conclusion (iii) by induction on $n$. Without loss of generality, we only show the case $i+2j+1=n$. With the help of Maple, from (2.7) and (2.9) and noting that the conclusions (i) and (ii), we obtain
\begin{eqnarray*}
\begin{cases}
I_{-1,9}=-\frac{768}{46189h^5}[(200h^3+3000h^2+2024h+512)I_{0,1}\\
\qquad\quad+(663h^4+326h^3+239h^2+64h)I_{1,1}],\\
I_{0,8}=-\frac{2048}{315h^5}(h+1)^4I_{2,0},\\
I_{1,7}=-\frac{64}{7293h^5}[(385h^3+1385h^2+1480h+512)I_{0,1}\\
\qquad\ \ +(139h^3+171h^2+64h)I_{1,1}],\\
I_{2,6}=-\frac{128}{35h^5}(h+1)^3I_{2,0},\\
I_{3,5}=-\frac{16}{3003h^5}[(480h^2+1000h+512)I_{0,1}+
(39h^3+111h^2+64h)I_{1,1}],\\
I_{4,4}=-\frac{32}{105h^5}(h+1)^2(h+8)I_{2,0},\\
I_{5,3}=-\frac{4}{1001h^5}[(77h^2+584h+512)I_{0,1}+(59h^2+64h)I_{1,1}],\\
I_{6,2}=\frac{4}{15h^5}(h+1)(3h+8)I_{2,0},\\
I_{7,1}=-\frac{1}{231h^5}[(232h+512)I_{0,1}+(15h^2+64h)I_{1,1}],\\
I_{8,0}=-\frac{1}{5h^5}(h^2+12h+16)I_{2,0},
\end{cases}
\end{eqnarray*}
which imply that the conclusion holds for $n=8$. Now assume that (iii) holds for $i+2l+1\leq k-1\, (k\geq9)$. For $i+2l+1=k$, if $k$ is an even number, then taking $(i,2l+1)=(-1,k+1)$ in (2.7) and $(i,2l+1)=(1,k-1),(3,k-3),\cdots,(k-3,3),(k-1,1)$ in (2.9), respectively, we have
\begin{eqnarray}
\mathbf{A}\left(\begin{matrix}
                I_{-1,k+1}\\
                 I_{1,k-1}\\
                 I_{3,k-3}\\
                 \vdots\\
                  I_{k-3,3}\\
                  I_{k-1,1}
                \end{matrix}\right)\ \
=\frac{1}{h}\left(\begin{matrix}
                \frac{4(k+1)}{5-3k}hI_{0,k-1}\\
               \frac{1}{3k-7}\big[(k-9)I_{-1,k-1}-4(k-4)I_{0,k-1}\big]\\
                \frac{1}{3k-9}\big[(k-7)I_{1,k-3}-4(k-4)I_{2,k-3}\big]\\
                                  \vdots\\
                  \frac{1}{2k-3}\big(2k-13)I_{k-5,3}-4(k-4)I_{k-4,3}\\
                \frac{1}{2k-5}\big(2k-11)I_{k-3,1}-4(k-4)I_{k-2,1}
                \end{matrix}\right),
\end{eqnarray}
where
\begin{eqnarray*}
\mathbf{A}=\left(\begin{matrix}
                1&\frac{4(k+1)}{5-3k}&0&\cdots&0&0\\
                 0&1&0&\cdots&0&0\\
                 0&0&1&\cdots&0&0\\
                 \vdots&\vdots&\vdots&\vdots&\vdots&\vdots\\
                  0&0&0&\cdots&1&0\\
                  0&0&0&\cdots&0&1
                \end{matrix}\right)\ \
\end{eqnarray*}
is a $\frac{k+2}{2}\times \frac{k+2}{2}$ matrix and $\det \mathbf{A}=1$.
If $k$ is an odd number taking $(i,2l+1)=(0,k)$ in (2.7) and $(i,2l+1)=(2,k-2),(4,k-4),\cdots,(k-3,3),(k-1,1)$ in (2.9), respectively, we have
\begin{eqnarray}
\mathbf{B}\left(\begin{matrix}
                I_{0,k}\\
                 I_{2,k-2}\\
                 I_{4,k-4}\\
                 \vdots\\
                  I_{k-3,3}\\
                  I_{k-1,1}
                \end{matrix}\right)\ \
=\frac{1}{h}\left(\begin{matrix}
                \frac{4k}{6-3k}hI_{1,k-2}\\
               \frac{1}{3k-8}\big[(k-8)I_{0,k-2}-4(k-4)I_{1,k-2}\big]\\
                \frac{1}{3k-10}\big[(k-6)I_{2,k-4}-4(k-4)I_{3,k-4}\big]\\
                                  \vdots\\
                  \frac{1}{2k-3}\big(2k-13)I_{k-5,3}-4(k-4)I_{k-4,3}\\
                \frac{1}{2k-5}\big(2k-11)I_{k-3,1}-4(k-4)I_{k-2,1}
                \end{matrix}\right),
\end{eqnarray}
where
\begin{eqnarray*}
\mathbf{B}=\left(\begin{matrix}
                1&\frac{4k}{6-3k}&0&\cdots&0&0\\
                 0&1&0&\cdots&0&0\\
                 0&0&1&\cdots&0&0\\
                 \vdots&\vdots&\vdots&\vdots&\vdots&\vdots\\
                  0&0&0&\cdots&1&0\\
                  0&0&0&\cdots&0&1
                \end{matrix}\right)\ \
\end{eqnarray*}
is a $\frac{k+1}{2}\times \frac{k+1}{2}$ matrix and $\det \mathbf{B}=1$. Hence, we can get that $I_{i,2l+1}$ can be expressed by $I_{0,1}$ and $I_{1,1}$ for $i+2l+1=k$ by the induction hypothesis.

From (2.15) and (2.16), we have for $(i,2l+1)=(-1,k+1)$ or $(i,2l+1)=(0,k)$
\begin{eqnarray*}
\begin{aligned}
I_{-1,k+1}(h)=&\frac{1}{h^{k-3}}\big[h\alpha^{(k-1)}(h)I_{0,1}+h\beta^{(k-1)}(h)I_{1,1}\big]\\
:=&\frac{1}{h^{k-3}}\big[\alpha^{(k)}(h)I_{0,1}+\beta^{(k)}(h)I_{1,1}\big],\quad  k \ \textup{even},\\
I_{0,k}(h)=&\frac{1}{h^{k-3}}\big[h\alpha^{(k-1)}(h)I_{0,1}+h\beta^{(k-1)}(h)I_{1,1}\big]\\
:=&\frac{1}{h^{k-3}}\big[\alpha^{(k)}(h)I_{0,1}+\beta^{(k)}(h)I_{1,1}\big], \, \quad k\ \textup{odd},\\
\end{aligned}
\end{eqnarray*}
where $\alpha^{(k-1)}(h)$ and $\beta^{(k-1)}(h)$ are polynomials in $h$. By the induction hypothesis we obtain that
$$\deg\alpha^{(k-1)}(h)\leq k-6,\ \ \deg\beta^{(k-1)}(h)\leq k-5.$$ Therefore,
$$\deg\alpha^{(k)}(h)\leq k-5,\ \ \deg\beta^{(k)}(h)\leq k-4.$$

In a similar way, we can prove the cases for $(i,2l+1)=(1,k-1),(3,k-3),\cdots,(k-3,3),(k-1,1)$ or $(i,2l+1)=(2,k-2),(4,k-4),\cdots,(k-3,3),(k-1,1)$.
This ends the proof.\quad $\lozenge$

\vskip 0.2 true cm

\noindent
{\bf Lemma 2.2.}\, {\it Suppose that $h\in(-1,0)$.}

 \vskip 0.2 true cm

\noindent
(i)\, If $n=2,3$, then
\begin{eqnarray}
M(h)=
{\alpha}(h)I_{0,1}(h)+{\beta}(h)I_{1,1}(h)+{\gamma}(h)I_{2,0}(h)
+{\delta}(h)I_{0,2}(h),
\end{eqnarray}
where
${\alpha}(h)$ is a constant, and ${\beta}(h)$, ${\gamma}(h)$ and ${\delta}(h)$ are polynomials of $h$ with $\deg {\beta}(h), \deg{\gamma}(h), \deg{\delta}(h)\leq1$.

 \vskip 0.2 true cm

\noindent
(ii)\, {If $4\leq n\leq7$, then
\begin{eqnarray*}
\begin{aligned}
M(h)=\frac{1}{h^{n-3}}[\alpha(h) I_{0,1}(h)+\beta(h) I_{1,1}(h)+\gamma(h) I_{2,0}(h)+\delta(h) I_{0,2}(h)],
\end{aligned}
\end{eqnarray*}
 where $\alpha(h)$, $\beta(h)$, $\gamma(h)$ and $\delta(h)$ are polynomials of $h$ with $\deg \alpha(h),\deg \delta(h)\leq n-4$ and $\deg \beta(h),\deg \gamma(h)\leq n-3$. }
  \vskip 0.2 true cm

\noindent
(iii)\, {If $n\geq8$, then
\begin{eqnarray*}
\begin{aligned}
M(h)=\frac{1}{h^{n-3}}[\alpha(h) I_{0,1}(h)+\beta(h) I_{1,1}(h)+\gamma(h) I_{2,0}(h)+\delta(h) I_{0,2}(h)],
\end{aligned}
\end{eqnarray*}
 where $\alpha(h)$, $\beta(h)$, $\gamma(h)$ and $\delta(h)$ are polynomials of $h$ with $\deg \alpha(h)\leq n-5$, $\deg \beta(h),\deg \gamma(h)\leq n-4$ and $\deg \delta(h)\leq 3$. }



 \vskip 0.2 true cm

\noindent
{\bf Lemma 2.3.}\, (i) {\it The vector function $(I_{0,1},I_{1,1})^T$ satisfies the following Picard-Fuchs equation
\begin{eqnarray}
\left(\begin{matrix}
                I_{0,1}\\
                 I_{1,1}
                \end{matrix}\right)
=\left(\begin{matrix}
                \frac{4}{5}h+\frac{16}{15}&\frac{4}{15}h\\
                 \frac{4}{3}&\frac{4}{3}h
                \end{matrix}\right)
                \left(\begin{matrix}
                I'_{0,1}\\
                 I'_{1,1}
                \end{matrix}\right).
\end{eqnarray}}
\noindent
(ii) {\it The vector function $(I_{2,0},I_{0,2})^T$ satisfies the following Picard-Fuchs equation
\begin{eqnarray}
\left(\begin{matrix}
                  I_{2,0}\\
                  I_{0,2}
                \end{matrix}\right)
=\left(\begin{matrix}
                2h+2&0\\
                 4h+4&h
                \end{matrix}\right)
                \left(\begin{matrix}
                  I'_{2,0}\\
                  I'_{0,2}
                \end{matrix}\right).
\end{eqnarray}}
 \vskip 0.2 true cm

\noindent
{\bf Proof.}\, From (1.3) we get
$$\frac{\partial y}{\partial h}=\frac{x^3}{y},$$
which implies
\begin{eqnarray}
I'_{i,j}=j\int_{\Gamma^+_h}x^{i-1}y^{j-2}dx.
\end{eqnarray}
Hence,
\begin{eqnarray}
I_{i,j}=\frac{1}{j+2}I'_{i-3,j+2}.
\end{eqnarray}
Multiplying both side of (2.20) by $h$, we have
\begin{eqnarray}
hI'_{i,j}=\frac{j}{2(j+2)}I'_{i-3,j+2}-2I'_{i-1,j}+I'_{i-2,j}.
\end{eqnarray}

From (2.3) and (2.20) we have for $j\geq1$
\begin{eqnarray}
\begin{aligned}
I_{i,j}=&\int_{\Gamma^+_h}x^{i-4}y^jdx=-\frac{j}{i-3}
\int_{\Gamma^+_h}x^{i-3}y^{j-1}dy\\
=&-\frac{j}{i-3}\int_{\Gamma^+_h}x^{i-3}y^{j-1}\frac{3hx^{2}+4x-1}{y}dx\\
=&-\frac{1}{i-3}(3hI'_{i,j}+4I'_{i-1,j}-I'_{i-2,j}).
\end{aligned}
\end{eqnarray}
(2.21)-(2.23) imply
\begin{eqnarray}
I_{i,j}=-\frac{4}{2i+j-6}\big(hI'_{i,j}+I'_{i-1,j}\big),\ j\geq1.
\end{eqnarray}

From (2.21) and noting that (2.14) we obtain
\begin{eqnarray*}
I_{2,0}=\frac{1}{2}I'_{-1,2}=\frac{2}{3}I_{2,0}+\frac{2}{3}(h+1)I'_{2,0}.
\end{eqnarray*}
Hence,
\begin{eqnarray}
I_{2,0}=2(h+1)I'_{2,0}.
\end{eqnarray}
From (2.24) we have
\begin{eqnarray}
\begin{aligned}
I_{0,1}=\frac{4}{5}(hI'_{0,1}+I'_{-1,1}),\ I_{1,1}=\frac{4}{3}(hI'_{1,1}+I'_{0,1}),\
I_{0,2}=hI'_{0,2}+I'_{-1,2},
\end{aligned}
\end{eqnarray}
and noting that (2.12) and (2.14) we obtain the conclusions (i) and (ii).
This ends the proof.\quad $\lozenge$

 \vskip 0.2 true cm

\noindent
{\bf Lemma 2.4.}\, {\it For $h\in(-1,0)$, $$I_{2,0}(h)=c_1\sqrt{h+1},\ \  I_{0,2}(h)=2c_1\sqrt{h+1}-c_1h\ln\frac{1-\sqrt{h+1}}{1+\sqrt{h+1}},$$
where $c_1$ is a nonzero constant.}
\vskip 0.2 true cm

\noindent
{\bf Proof.}  From (2.19) we have $I_{2,0}(h)=c_1\sqrt{h+1}$, where $c_1$ is a constant. Therefore, we have for $h\in(-1,0)$
$$
I_{0,2}(h)=c_2h+2c_1\sqrt{h+1}-c_1h\ln\frac{1-\sqrt{h+1}}{1+\sqrt{h+1}}
$$
where $c_2$ is a constant. Since $I_{0,2}(-1)=0$, we have $c_2=0$.  Hence, $I_{0,2}(h)=2c_1\sqrt{h+1}-c_1h\ln\frac{1-\sqrt{h+1}}{1+\sqrt{h+1}}$. This ends the proof.\quad $\lozenge$
\vskip 0.2 true cm

Taking $(i,j)=(4,1),(3,1)$ in (2.9) respectively and bearing in mind that (2.2), we get
$$\begin{aligned}
I_{3,1}(h)=-\frac{1}{h}I_{1,1}(h),\ \ I_{4,1}(h)=-\frac{1}{5h}[I_{2,1}(h)+4I_{3,1}(h)].
\end{aligned}$$
Hence, $I_{0,1}(h)=h^2I_{4,1}(h)$. Using Green formula, we have $$I_{4,1}(h)=\int_{\Gamma^+_h}ydx=\oint_{\Gamma^+_h\cup\overrightarrow{AB}}ydx=\iint\limits_Ddxdy\neq0,$$
where $D$ is the interior of $\Gamma_{h}^+\cup \overrightarrow{AB}$, see Fig.\,1. Thus, $I_{0,1}(h)\neq0$ for $h\in(-1,0)$. Noting that $\frac{\partial y}{\partial h}=x^3y^{-1}$ and $dx=xydt$, we have
$$I'_{0,1}(h)=\int_{\Gamma^+_h}x^{-4}\frac{\partial y}{\partial h}dx=
\int_0^{t_0}dt\neq0,$$
where $t_0$ is the time from the left end point to right end point of $\Gamma^+_h$. So we can get the following lemma.

 \vskip 0.2 true cm

\noindent
{\bf Lemma 2.5.} {\it Let $\omega_1(h)=\frac{I_{1,1}(h)}{I_{0,1}(h)}$ and $\omega_2(h)=\frac{I'_{1,1}(h)}{I'_{0,1}(h)}$ for $h\in (-1,0)$, then $\omega_1(h)$ and $\omega_2(h)$ satisfy the following Riccati equations
\begin{eqnarray}
G(h)\omega'_1(h)=\frac{1}{4}h\omega^2_1(h)-\frac{1}{2}(h-2)
\omega_1(h)-\frac{5}{4}
\end{eqnarray}
and
\begin{eqnarray}
G(h)\omega'_2(h)=-\frac{1}{4}h\omega^2_2(h)-\frac{1}{2}h
\omega_2(h)-\frac{1}{4},
\end{eqnarray}
respectively, where $G(h)=h(h+1)$.}

 \vskip 0.2 true cm

\noindent
{\bf Proof.}\, From (2.18), we have
\begin{eqnarray*}G(h)\left(\begin{matrix}
                  I'_{0,1}(h)\\
                  I'_{1,1}(h)
                \end{matrix}\right)
=\left(\begin{matrix}
                \frac{5}{4}h&-\frac{1}{4}h\\
                -\frac{5}{4}&\frac{3}{4}h+1
                \end{matrix}\right)
                \left(\begin{matrix}
                  I_{0,1}(h)\\
                  I_{1,1}(h)
                \end{matrix}\right)
\end{eqnarray*}
and
\begin{eqnarray*}G(h)\left(\begin{matrix}
                  I''_{0,1}(h)\\
                  I''_{1,1}(h)
                \end{matrix}\right)
=\left(\begin{matrix}
                \frac{1}{4}h&-\frac{1}{4}h\\
                -\frac{1}{4}&-\frac{1}{4}h
                \end{matrix}\right)
                \left(\begin{matrix}
                  I'_{0,1}(h)\\
                  I'_{1,1}(h)
                \end{matrix}\right),
\end{eqnarray*}
where $G(h)=h(h+1)$. Noting that $G(h)\neq0$ for $h\in(-1,0)$ and $$\omega'_1(h)=\frac{I'_{1,1}(h)}{I_{0,1}(h)}-
\omega_1(h)\frac{I'_{0,1}(h)}{I_{0,1}(h)},\ \omega'_2(h)=\frac{I''_{1,1}(h)}{I'_{0,1}(h)}-
\omega_2(h)\frac{I''_{0,1}(h)}{I'_{0,1}(h)},$$
 we obtain (2.27) and (2.28). This ends the proof. \quad $\lozenge$

\section{Proof of the main results}
 \setcounter{equation}{0}
\renewcommand\theequation{3.\arabic{equation}}

In order to prove the Theorem 1.1, we first introduce some helpful results in the literature. Let $V$ be a finite-dimensional vector space of functions, real-analytic on an open interval $\mathbb{I}$.

\vskip 0.2 true cm

\noindent
{\bf Definition 3.1}\, {\it \cite{GI}. We say that $S$ is a Chebyshev space, provided that each non-zero function in $S$ has at most $\dim(S)-1$ zeros, counted with multiplicity.}

 \vskip 0.2 true cm
\noindent
{\bf Proposition 3.1}\, {\it \cite{GI}. The solution space $S$ of a second order linear analytic differential equation $$x''+a_1(t)x'+a_2(t)x=0$$
on an open interval $\mathbb{I}$ is a Chebyshev space if and only if there exists a nowhere vanishing solution $x_0(t)\in S\ (x_0(t)\neq0,\ \forall t\in\mathbb{I}$).}

\vskip 0.2 true cm
\noindent
{\bf Proposition 3.2}\, {\it \cite{GI}. Suppose the solution space of the homogeneous equation $$x''+a_1(t)x'+a_2(t)x=0$$ is a Chebyshev space and let $R(t)$ be an analytic function on $\mathbb{I}$ having $l$ zeros (counted with multiplicity). Then every solution $x(t)$ of the non-homogeneous equation
$$x''+a_1(t)x'+a_2(t)x=R(t)$$
has at most $l+2$ zeros on $\mathbb{I}$.}
\vskip 0.2 true cm

In the following we denote by $\#\{\varphi(h)=0,h\in(a,b)\}$ the number of isolated zeros of $\varphi(h)$ on $(a,b)$ taking into account the multiplicity, and we also denote by $\Theta_k(h)$ the polynomial of degree at most $k$.
 \vskip 0.2 true cm

\noindent
{\bf Lemma 3.1.}\,{\it Suppose that $h\in(-1,0)$.}
 \vskip 0.2 true cm

\noindent
 (i)\, {\it If $n=2,3$, then there exist polynomials $P^1_2(h)$, $P^1_1(h)$ and $P^1_0(h)$ of $h$ with degree respectively $4$, $3$ and $2$ such that $L^1(h)\Phi(h)=0$.}
  \vskip 0.2 true cm

\noindent
  (ii)\, {\it If $4\leq n\leq7$, then there exist polynomials $P^2_2(h)$, $P^2_1(h)$ and $P^2_0(h)$ of $h$ with degree respectively $2n-4$, $2n-5$ and $2n-6$ such that $L^2(h)\Phi(h)=0$.}

\vskip 0.2 true cm

\noindent
  (iii)\, {\it If $n\geq8$, then there exist polynomials $P^3_2(h)$, $P^3_1(h)$ and $P^3_0(h)$ of $h$ with degree respectively $2n-6$, $2n-7$ and $2n-8$ such that $L^3(h)\Phi(h)=0$, where $$\Phi(h)=\alpha(h)I_{0,1}(h)+\beta(h)I_{1,1}(h),$$
\begin{eqnarray}
L^i(h)=P^i_2(h)\frac{d^2}{dh^2}+P^i_1(h)\frac{d}{dh}+P^i_0(h),\ i=1,2,3.
\end{eqnarray}}

\noindent
{\bf Proof.}\, Without loss of generality, we only prove (iii). (i) and (ii) can be shown similarly. By (2.18), we have
$$V'(h)=(E-B)^{-1}(Bh+C)V''(h),$$
where $V(h)=(I_{0,1}(h),I_{1,1}(h))^T$, and
$$E=\left(\begin{matrix}
                1&0\\
                 0&1
                \end{matrix}\right),\ \
B=\left(\begin{matrix}
                \frac{4}{5}&\frac{4}{15}\\
                 0&\frac{4}{3}
                \end{matrix}\right),\ \
C=\left(\begin{matrix}
                \frac{16}{15}&0\\
                 \frac{4}{3}&0
                \end{matrix}\right).$$
Hence,
\begin{eqnarray*}
\begin{aligned}
\Phi(h)=&\tau(h) V(h)=\tau(h)(Bh+C)V'(h)\\
=&\tau(h)(Bh+C)(E-B)^{-1}(Bh+C) V''(h)\\
:=&\Theta_{n-3}(h)I''_{0,1}(h)+\Theta_{n-2}(h)I''_{1,1}(h),
\end{aligned}
\end{eqnarray*}
where $\tau(h)=(\alpha(h),\beta(h))$, $\Theta_{n-3}(h)$ denotes a polynomial in $h$ of degree at most $n-3$ and etc..
For $\Phi'(h)$, we have
\begin{eqnarray*}
\begin{aligned}
\Phi'(h)=&\tau'(h)V(h)+\tau(h)V'(h)\\
=&\Big[\tau'(h)(Bh+C)+\tau(h)\Big](E-B)^{-1}(Bh+C) V''(h)\\
:=&\Theta_{n-4}(h)I''_{0,1}(h)+\Theta_{n-3}(h)I''_{1,1}(h).
\end{aligned}
\end{eqnarray*}
In a similar way, we have
\begin{eqnarray*}
\begin{aligned}
\Phi''(h)=\Theta_{n-5}(h)I''_{0,1}(h)+\Theta_{n-4}(h)I''_{1,1}(h).
\end{aligned}
\end{eqnarray*}

Next, suppose that
\begin{eqnarray}
P_2(h)=\sum\limits_{k=0}^{2n-6}p_{2,k}h^k,\ \ P_1(h)=\sum\limits_{m=0}^{2n-7}p_{1,m}h^m, \ \ P_0(h)=\sum\limits_{l=0}^{2n-8}p_{0,l}h^l
\end{eqnarray}
are polynomials of $h$ with coefficients $p_{2,k},\ p_{1,m}$ and $p_{0,l}$ to be determined such that $L(h)\Phi(h)=0$ for
\begin{eqnarray}
0\leq k\leq 2n-6,\ \  0\leq m\leq 2n-7,\ \ 0\leq l\leq 2n-8.
\end{eqnarray}
By straightforward computation, we have
\begin{eqnarray*}
\begin{aligned}
L(h)\Phi(h)=&P_2(h)\Phi''(h)+P_1(h)\Phi'(h)+P_0(h)\Phi(h)\\[0.2cm]
=&\Big[P_2(h)\Theta_{n-5}(h)+P_1(h)\Theta_{n-4}(h)+P_0(h)\Theta_{n-3}(h)\Big]I''_{0,1}(h)\\
&~+\Big[P_2(h)\Theta_{n-4}(h)+P_1(h)\Theta_{n-3}(h)+P_0(h)\Theta_{n-2}(h)\Big]I''_{1,1}(h)\\
:=&X(h)I''_{0,1}(h)+Y(h)I''_{1,1}(h),
\end{aligned}
\end{eqnarray*}
where $X(h)$ and $Y(h)$ are polynomials of $h$ with $\deg X(h)\leq 3n-11$ and $\deg Y(h)\leq 3n-10$.

 Let
$$X(h)=\sum\limits_{i=0}^{3n-11}x_ih^i,\ \ Y(h)=\sum\limits_{j=0}^{3n-10}y_jh^j,$$
where $x_i$ and $y_j$ are expressed by $p_{2,k}$, $p_{1,m}$ and $p_{0,l}$ of (3.2) linearly, $k$, $m$ and $l$ satisfy (3.3). Let
\begin{eqnarray}
x_i=0,\ \ y_j=0, \ \ 0\leq i\leq 3n-11,\ \ 0\leq j\leq 3n-10,
\end{eqnarray}
then system (3.4) is a homogenous linear equations with $6n-19$ equations about $6n-18$ variables of $p_{2,k}$, $p_{1,m}$ and $p_{0,l}$ for $k$, $m$ and $l$ satisfy (3.3). It follows that from the theory of linear algebra that there exist $p_{2,k}$, $p_{1,m}$ and $p_{0,l}$ such that (3.4) holds, which yields $L(h)\Phi(h)=0$. This ends the proof.\quad $\lozenge$
\vskip 0.2 true cm

 \noindent
{\bf Lemma 3.2.}\, {\it Let $\Phi(h)=\alpha(h)I_{0,1}(h)+\beta(h)I_{1,1}(h)$.}
\vskip 0.2 true cm

\noindent
(i)\, {\it If $n=2,3$, then $\Phi(h)$ has at most $4$ zeros on $(-1,0)$, taking into account the multiplicity.}
\vskip 0.2 true cm

\noindent
(ii)\, {\it If $4\leq n\leq7$, then $\Phi(h)$ has at most $3n-8$ zeros on $(-1,0)$, taking into account the multiplicity.}
\vskip 0.2 true cm

\noindent
(iii)\, {\it If $n\geq8$, then $\Phi(h)$ has at most $3n-11$ zeros on $(-1,0)$, taking into account the multiplicity.}
\vskip 0.2 true cm

\noindent
{\bf Proof.}\, We only prove (iii). (i) and (ii) can be proved in a similar way. Let $\chi_1(h)=\alpha(h)+\beta(h)\omega_1(h)$, so $\Phi(h)=I_{0,1}(h)\chi_1(h)$ which implies
$$\#\{\Phi(h)=0,h\in(-1,0)\}=\#\{\chi_1(h)=0,h\in(-1,0)\}.$$
By (2.27) we know that $\chi_1(h)$ satisfies
\begin{eqnarray}
G(h)\beta(h)\chi'_1(h)=\frac{1}{4}h\chi_1(h)^2+F_1(h)\chi_1(h)+F_0(h)
\end{eqnarray}
 with $\deg F_0(h)\leq2n-8$. Recall that the inequality (4.8) in \cite{ZZ} is
$$\nu\leq\sigma+\lambda+1,$$
where $\nu$, $\sigma$ and  $\lambda$ correspond here to $\#\{\chi_1(h)=0,h\in(-1,0)\}$, $\#\{F_0(h)=0,h\in(-1,0)\}$ and $\#\{\beta(h)=0,h\in(-1,0)\}$, respectively. Hence, we have for $h\in(-1,0)$
$$\#\{\chi_1(h)=0\}\leq\#\{\beta(h)=0\}+\#\{F_0(h)=0\}+1\leq3n-11.$$
Hence,
$$\#\{\Phi(h)=0,h\in(-1,0)\}=\#\{\chi_1(h)=0,h\in(-1,0)\}\leq3n-11.$$
This completes the proof.\quad $\lozenge$
\vskip 0.2 true cm

\noindent
{\bf Proof of the Theorem 1.1.} We only prove (iii). (i) and (ii) can be proved similarly.

 Let $M_1(h)=h^{n-3}M(h)$, then $M_1(h)$ has the same zeros as $M(h)$ on $(-1,0)$. For the sake of clearness, we split the proof into three steps.
\vskip 0.2 true cm
 (1) For $h\in(-1,0)$, $L^3(h)M_1(h)=R(h)$, where $L^3(h)$ is defined by (3.1),
 \begin{eqnarray}
 R(h)=\Theta_{2n-4}(h)\ln\frac{1-\sqrt{h+1}}{1+\sqrt{h+1}}+\Theta_{3n-9}(h)\frac{1}{h(h+1)^\frac{3}{2}}.
 \end{eqnarray}
\vskip 0.2 true cm

In fact, from Lemma 2.4, we have
\begin{eqnarray}
\begin{aligned}
\Psi(h):=&\gamma(h)I_{2,0}(h)+\delta(h)I_{0,2}(h)\\
=&c_1[\gamma(h)+2\delta(h)]\sqrt{h+1}-c_1h\delta(h)
\ln\frac{1-\sqrt{h+1}}{1+\sqrt{h+1}}\\
:=&\Theta_{n-4}(h)\sqrt{h+1}+h\Theta_3(h)\ln\frac{1-\sqrt{h+1}}{1+\sqrt{h+1}},\\
\Psi'(h)=&\Theta_{n-4}(h)\frac{1}{\sqrt{h+1}}+\Theta_3(h)\ln\frac{1-\sqrt{h+1}}{1+\sqrt{h+1}},\\
\Psi''(h)=&\Theta_{n-3}(h)\frac{1}{h(h+1)^\frac{3}{2}}+\Theta_2(h)\ln\frac{1-\sqrt{h+1}}{1+\sqrt{h+1}}.
\end{aligned}
\end{eqnarray}
 From Lemma 3.1 (iii), we have
\begin{eqnarray}
L^3(h)M_1(h)=L^3(h)\Psi(h)=P^3_2(h)\Psi''(h)+P^3_1(h)\Psi'(h)+P^3_0(h)\Psi(h).
\end{eqnarray}
Substituting (3.7) into (3.8) gives (3.6).
\vskip 0.2 true cm

 (2) Zeros of $R(h)$ for $h\in(-1,0)$.
\vskip 0.2 true cm

Denote that $U=\{h\in(-1,0)|\Theta_{2n-4}(h)=0\}$. For $h\in (-1,0)\backslash U$, by detailed computations, we get
\begin{eqnarray}
\begin{aligned}
\Big(\frac{R(h)}{\Theta_{2n-4}(h)}\Big)'=\frac{\Theta_{5n-12}(h)}{\Theta_{2n-4}^2(h)h^2(h+1)^\frac{5}{2}}.
\end{aligned}
\end{eqnarray}
Since $h^2(h+1)^\frac{5}{2}\neq 0$ for $h\in(-1,0)$, we have
\begin{eqnarray}
\#\{R(h)=0,h\in(-1,0)\}\leq 7n-15.
\end{eqnarray}

\vskip 0.1 true cm

(3) Zeros of $M(h)$ for $h\in(-1,0)$.

\vskip 0.2 true cm

By Lemma 3.2, we have $\Phi(h)$ has at most $3n-11$ zeros on $(-1,0)$.
We assume that
$$P^3_2(\tilde{h}_i)=0,\ \Phi(\bar{h}_j)=0,\ \tilde{h}_i, \bar{h}_j\in(-1,0),\ 1\leq i\leq 2n-6,\ 1\leq j\leq3n-11.$$
Denote $\tilde{h}_i$ and $\bar{h}_j$ as $h_m^*$, and reorder them such that $h_m^*<h_{m+1}^*$ for $m=1,2,\cdots,5n-17$. Let
$$\Delta_s=(h_s^*,h_{s+1}^*),\ s=0,1,\cdots,5n-17,$$
where $h_0^*=-1$, $h_{5n-16}^*=0$. Then $P^3_2(h)\neq0$ and $\Phi(h)\neq0$ for $h\in\Delta_s$ and $L^3(h)\Phi(h)=0$. By Proposition 3.1, the solution space of
$$L^3(h)=P^3_2(h)\Big(\frac{d^2}{dh^2}+\frac{P^3_1(h)}{P_2(h)}\frac{d}{dh}+\frac{P^3_0(h)}{P_2(h)}\Big)=0$$
is a Chebyshev space on $\Delta_s$. By Proposition 3.2, $M_1(h)$ has at most $2+l_s$ zeros for $h\in\Delta_s$, where $l_s$ is the number of zeros of $R(h)$ on $\Delta_s$. Therefore, we obtain for $h\in(-1,0)$
\begin{eqnarray*}
\begin{aligned}
\#\{M(h)=0\}=&\#\{M_1(h)=0\}\\
\leq&\#\{R(h)=0\}+2\cdot \textup{the number of the intervals of } \Delta_s\\
&+\textup{the number of the end points of }\Delta_s \\
\leq&22n-64.
\end{aligned}
\end{eqnarray*}
\vskip 0.2 true cm

\noindent
{\bf Proof of the Theorem 1.2.}
If $a^+_{i,j}=a^-_{i,j}$ and $b^+_{i,j}=b^-_{i,j}$, that is, the systems (1.4) is smooth. Since $\Gamma_h$ is symmetric with respect to $x$-axis for $h\in(-1,0)$, $A_{i,2l}(h)=\oint_{\Gamma_h}x^{i-4}y^{2l}dx=0$, $l=0,1,2,\cdots$, where
$$\Gamma_h=\Gamma^+_h\cup\Gamma^-_h,\ \ A_{i,j}(h)=I_{i,j}(h)+J_{i,j}(h).$$
Hence, from Lemma 2.2 we have
\begin{eqnarray*}
M(h)=\begin{cases}
\frac{1}{h^{n-3}}[\tilde{\alpha}(h) A_{0,1}(h)+\tilde{\beta}(h) A_{1,1}(h)],\ n=2,3,\\
\frac{1}{h^{n-3}}[\alpha(h) A_{0,1}(h)+{\beta}(h) A_{1,1}(h)],\ n\geq4,
\end{cases}
\end{eqnarray*}
 where $\tilde{\alpha}(h)$ is a constant, and $\tilde{\beta}(h)$, $\alpha(h)$ and $\beta(h)$ are polynomials of $h$ with $\deg \tilde{\beta}(h)\leq 1$, $\deg \alpha(h)\leq n-4$ and $\deg \beta(h)\leq n-3$.  By the same proof of
Lemma 3.2, we have
 $$\#\{M(h)=0,h\in(-1,0)\}\leq\begin{cases}4,\qquad\, \ n=2,3,\\
 3n-8,\ n\geq4.\end{cases}$$

 \vskip 0.2 true cm

\noindent
{\bf Acknowledgment}
 \vskip 0.2 true cm

\noindent
Supported by Higher Educational Science Program of Ningxia(NGY201789), National Natural Science Foundation of China(11701306,11671040,11601250), Construction of First-class Disciplines of Higher Education of Ningxia(pedagogy)(NXYLXK2017B11), Key Program of Ningxia Normal University(NXSFZD1708) and Science and Technology Pillar Program of Ningxia(KJ[2015]26(4)).

\end{document}